\documentclass{amsproc}

\usepackage{amssymb}

\newtheorem{theorem}{Theorem}[section]

\newtheorem{proposition}[theorem]{Proposition}

\theoremstyle{definition}
\newtheorem{definition}[theorem]{Definition}
\newtheorem{example}[theorem]{Example}

\theoremstyle{remark}
\newtheorem{remark}[theorem]{Remark}

\numberwithin{equation}{section}

 \DeclareMathOperator{\comp}{\#}
 
 \newcommand{\N}{\mathbf{N}}
 
 \newcommand{\Z}{\mathbf{Z}}

\begin{document}

\title{Simple omega-categories and chain complexes}

\author{Richard Steiner}

\address{Department of Mathematics\\University of Glasgow\\
University Gardens\\Glasgow\\Scotland G12 8QW}

\email{r.steiner@maths.gla.ac.uk}



\begin{abstract}
The category of strict omega-categories has an important full
subcategory whose objects are the simple omega-categories freely
generated by planar trees or by globular cardinals. We give a
simple description of this subcategory in terms of chain
complexes, and we give a similar description of the opposite
category, the category of finite discs, in terms of cochain
complexes. Berger has shown that the category of simple
omega-categories has a filtration by iterated wreath products of
the simplex category. We generalise his result by considering
wreath products of categories of chain complexes over the simplex
category.
\end{abstract}

\maketitle

\section{Introduction} \label{S1}

The category of strict $\omega$-categories has an important full
subcategory~$\Theta$, whose objects are the $\omega$-categories
freely generated by planar trees in the sense of Batanin \cite{Ba}
(see also~\cite{BS}). One can regard~$\Theta$ as the theory of
strict $\omega$-categories in the sense of universal algebra, and
it has been used in the study of weak $\omega$-categories by
Batanin~\cite{Ba} and Joyal~\cite{J}. The subcategories~$\Theta_n$
of~$\Theta$ consisting of $n$-categories for finite values of~$n$
have been applied to iterated loop spaces by Berger~\cite{Be}.

In this paper, the objects of~$\Theta$ will be called simple
$\omega$-categories, following Makkai and Zawadowski~\cite{MZ};
they were called Batanin cells in~\cite{J}, which is also the
source of the notation~$\Theta$. The generating structures for the
objects of~$\Theta$ were called globular cardinals by
Street~\cite{StrPet}. There are several ways to describe the
category~$\Theta$, most of which give rather complicated
descriptions of the morphisms; the main object of this paper is to
give a simple description of the objects and morphisms of~$\Theta$
in terms of chain complexes and chain maps. We give a new
description of the generating structures for the objects
of~$\Theta$ in Section~\ref{S2} and the description of the
category~$\Theta$ itself in Section~\ref{S3}. The method gives a
similar description for the morphisms between simple
$\omega$-categories and certain other $\omega$-categories, for
example the oriented simplexes or orientals of
Street~\cite{StrAlg} (see also~\cite{SteOr}). The opposite
category to~$\Theta$ has been studied under the name of the
category of finite discs \cite{J}, \cite{MZ}; in Section~\ref{S4}
we obtain a simple description of this category in terms of
cochain complexes and cochain maps. Berger~\cite{Be} has shown
that the subcategories~$\Theta_n$ of~$\Theta$ are equivalent to
iterated wreath products of the simplex category; in
Section~\ref{S5} we generalise his result by studying wreath
products of certain categories of chain complexes over the simplex
category.

\section{Simple omega-categories and graded ordered sets} \label{S2}

In this section we give the definition of simple
$\omega$-categories, and we show that they are generated by
certain graded ordered sets. The term simple $\omega$-category is
due to Makkai and Zawadowski \cite{MZ}, but our treatment is based
on the work of Street \cite{StrPet}.

All of the $\omega$-categories in this paper are strict
$\omega$-categories. We regard an $\omega$-category as a single
set with an infinite sequence of partially defined binary
composition operations, each of which makes it the set of
morphisms in a small category. The composition operations are
denoted $\comp_0,\comp_1,\ldots\,$, the left identity of an
element~$x$ under~$\comp_n$ is denoted $d_n^- x$, and the right
identity of~$x$ under~$\comp_n$ is denoted $d_n^+ x$. The category
structures commute with one another, with the special feature that
\[
 d_m^\beta d_n^\alpha x=d_n^\alpha d_m^\beta x=d_m^\beta x\
 \text{for $m<n$};
\]
in other words, if $x$~is an identity for some~$\comp_m$, then it
is also an identity for $\comp_{m+1},\comp_{m+2},\ldots\,$. There
is a final axiom saying that every element~$x$ is an identity for
some operation~$\comp_m$.

Simple $\omega$-categories are the free $\omega$-categories on a
particular kind of globular set, called a simple globular set or a
globular cardinal. We will use the following notation and
definitions.

A \emph{globular set} (sometimes called an \emph{$\omega$-graph})
is a set in which: every element~$x$ is assigned a nonnegative
integer dimension~$|x|$; every positive-dimensional element~$|x|$
is assigned two $(|x|-1)$-dimensional elements called its source
and target and here denoted $\partial^- x$ and $\partial^+ x$; if
$|x|\geq 2$ then
\[
 \partial^-\partial^- x=\partial^-\partial^+ x,\
 \partial^+\partial^- x=\partial^+\partial^+ x.
\]
The \emph{free $\omega$-category on a globular set} is the
$\omega$-category with the following presentation: the generators
are the members of the globular set; if $x$~is an $n$-dimensional
generator then there are relations $d_n^- x=d_n^+ x=x$; if $x$~is
an $n$-dimensional generator with $n>0$ then there are relations
$d_{n-1}^- x=\partial^- x$ and $d_{n-1}^+ x=\partial^+ x$. A
\emph{simple globular set} or \emph{globular cardinal} is a
non-empty finite globular set such that the transitive closure of
the relation given by $\partial^- x<x$ and $x<\partial^+ x$ is a
total ordering. A \emph{simple $\omega$-category} is the free
$\omega$-category on a simple globular set.

Since a simple globular set is finite, the transitive closure
condition yields the following result.

\begin{proposition} \label{2.1}
Let $x$~and~$y$ be consecutive elements in the ordering of a
simple globular set with $x<y$. Then $x=\partial^- y$ or
$\partial^+ x=y$. The dimensions of $x$~and~$y$ differ by~$1$.
\end{proposition}

We deduce that the globular structure can be recovered from the
ordering and the dimension function as follows.

\begin{proposition} \label{2.2}
Let $x$ be a positive-dimensional element in a simple globular
set. Then $\partial^- x$ is the last $(|x|-1)$-dimensional element
before~$x$ in the ordering, and $\partial^+ x$ is the first
$(|x|-1)$-dimensional element after~$x$ in the ordering.
\end{proposition}

\begin{proof}
We will prove the result for $\partial^- x$. Let the dimension
of~$x$ be $n+1$. Certainly $\partial^- x$ is an $n$-dimensional
element coming before~$x$ in the ordering, so there really is a
last $n$-dimensional element~$a$ before~$x$, and we must show that
$a=\partial^- x$.

Let $b$ be the immediate successor of~$a$. From
Proposition~\ref{2.1}, the dimensions of the elements~$y$ such
that $b\leq y\leq x$ form a consecutive set of integers. This set
contains $n+1$ and does not contain~$n$; therefore $|y|>n$ for
$b\leq y\leq x$. In particular, Proposition~\ref{2.1} now shows
that $|b|=n+1$ and $a=\partial^- b$.

For consecutive elements $y$~and~$z$ with $b\leq y<z\leq x$ we
have $|y|>n$ and $|z|>n$, and we claim further that
$(\partial^-)^{|y|-n}y=(\partial^-)^{|z|-n}z$; indeed this is
trivially true if $y=\partial^- z$, and it follows from the
identity $\partial^-\partial^-=\partial^-\partial^+$ if
$\partial^+y=z$.

We now see that $y\mapsto(\partial^-)^{|y|-n}y$ is constant on the
interval $b\leq y\leq x$. In particular
$(\partial^-)^{|b|-n}b=(\partial^-)^{|x|-n}x$. Since $a=\partial^-
b$ and $|b|=|x|=n+1$, this gives us the required equality
$a=\partial^- x$.
\end{proof}

It follows from Proposition~\ref{2.2} that simple globular sets
are equivalent to ordered sets with suitable dimension functions;
in other words they are equivalent to certain graded ordered sets.
The graded ordered sets that can occur turn out to be the
continuously graded ones in the sense of the following definition.

\begin{definition} \label{2.3}
A \emph{continuously graded ordered set} is a non-empty finite
ordered set, together with a function assigning a nonnegative
integer dimension~$|x|$ to each element~$x$, such that the first
and last elements have dimension zero and such that consecutive
elements have dimensions differing by~$1$.
\end{definition}

For example, there is a continuously graded ordered set such that
the dimensions of its elements in order are
\[
 0,1,2,1,2,3,4,3,2,3,4,3,2,1,0,1,0.
\]

The main result is now as follows.

\begin{theorem} \label{2.4}
Let $X$ be a continuously graded ordered set. Then there are
well-defined functions $\partial^-$~and~$\partial^+$ on the
positive-dimensional elements of~$X$ such that $\partial^- x$ is
the last $(|x|-1)$-dimensional element preceding~$x$ and
$\partial^+ x$ is the first $(|x|-1)$-dimensional element
following~$x$, and these functions make~$X$ into a simple globular
set. Every simple globular set arises from a continuously graded
ordered set in this way.
\end{theorem}

\begin{proof}
By the conditions in Definition~\ref{2.3}, if $x$~is a
positive-dimensional element of~$X$ then there is at least one
$(|x|-1)$-dimensional element before~$x$ and at least one
$(|x|-1)$-dimensional element after~$x$. Therefore the functions
$\partial^-$~and~$\partial^+$ are well-defined. If $x$~is at least
$2$-dimensional, then the conditions of Definition~\ref{2.3} imply
that  $|y|\geq |x|-1$ for $\partial^- x\leq y\leq\partial^+ x$,
and it follows that $\partial^\alpha\partial^-
x=\partial^\alpha\partial^+ x$ for each sign~$\alpha$; therefore
$X$~is a globular set. By construction, $\partial^- x<x<\partial^+
x$ for all positive-dimensional~$x$; also, if $x$~and~$y$ are
consecutive elements with $x<y$, then $|x|=|y|-1$ or $|y|=|x|-1$,
so $x=\partial^- y$ or $\partial^+ x=y$. Therefore the ordering is
the transitive closure obtained from the functions
$\partial^-$~and~$\partial^+$ in the way required for a simple
globular set, and it follows that $X$~is indeed a simple globular
set.

Conversely, let $X$ be a simple globular set. Since $\partial^-
x<x<\partial^+ x$ whenever $x$~has positive dimension, it follows
that the first and last elements are zero-dimensional. Also, by
Proposition~\ref{2.1}, consecutive elements have dimensions
differing by~$1$, so $X$~is a continuously graded ordered set, and
it follows from Proposition~\ref{2.2} that the functions
$\partial^-$~and~$\partial^+$ come from the ordering and
dimensions in the way described.

This completes the proof.
\end{proof}

Since simple globular sets are equivalent to continuously graded
ordered sets, their isomorphism classes can be indexed by the
non-empty finite sequences of nonnegative integers beginning and
ending with~$0$ and with adjacent terms differing by~$1$. The
isomorphism classes of simple $\omega$-categories can therefore be
indexed by these sequences as well. Without loss of information, a
sequence of this kind can be replaced by the subsequence
consisting of its maxima (terms equal to~$n$ and not adjacent to
$n+1$) and its internal minima (terms equal to~$n$ adjacent to
$n+1$ on both sides); for example the sequence
\[
 (0,1,2,1,2,3,4,3,2,3,4,3,2,1,0,1,0)
\]
can be replaced by
\[
 (2,1,4,2,4,0,1).
\]
The original sequence can be recovered from the subsequence by
interpolation. Note that the sequence~$(0)$ yields the
subsequence~$(0)$, but the initial and final zeros are omitted
from the subsequence in all other cases. The subsequences that
occur are the non-empty finite sequences of nonnegative integers
\[
 (u_0,v_1,u_1,v_2,u_2,\ldots,u_{k-1},v_k,u_k)
\]
such that
\[
 u_0>v_1,\ v_1<u_1,\ u_1>v_2,\ v_2<u_2,\ \ldots,\
 u_{k-1}>v_k,\ v_k<u_k;
\]
in other words they are the up-and-down vectors used to index
simple $\omega$-categories by Makkai and Zawadowski (\cite{MZ},
2.3). They may also be used to index planar trees with a
distinguished vertex and a distinguished maximal path starting at
that vertex. Indeed, given a tree with this structure, let
$P_0,\ldots,P_k$ be the maximal paths starting at the
distinguished vertex listed in clockwise order starting with the
distinguished path, let $u_i$ be the number of edges in~$P_i$, and
let $v_i$ be the number of edges in $P_{i-1}\cap P_i$; then
$(u_0,v_1,u_1,\ldots,u_k)$ is an up-and-down vector and every
up-and-down vector comes from a planar tree in this way. This
explains the indexing of simple $\omega$-categories by planar
trees (\cite{Ba}, \cite{BS}, \cite{StrPet}). Berger~\cite{Be}
calls these structures level-trees. In more purely combinatorial
terms, they are abstract trees with a distinguished root vertex
and a total ordering at each vertex on the set of edges at that
vertex pointing away from the root. To summarise, we have the
following result.

\begin{theorem} \label{2.5}
There are explicit bijections between the following sets\textup{:}
the set of isomorphism classes of simple globular sets\textup{;}
the set of isomorphism classes of continuously graded ordered
sets\textup{;} the set of up-and-down vectors\textup{;} the set of
isomorphism classes of level-trees.
\end{theorem}

\section{Simple augmented directed complexes} \label{S3}

We will now recall from \cite{SteOm} the theory associating
$\omega$-categories and chain complexes, and we will use it to
describe the category of simple $\omega$-categories in terms of
chain complexes. All our chain complexes will be augmented chain
complexes of abelian groups concentrated in nonnegative
dimensions. We recall that an \emph{augmented directed complex} is
a chain complex~$K$ of this type together with a distinguished
submonoid~$K_q^*$ of~$K_q$ for each chain group~$K_q$. A morphism
of augmented directed complexes from~$K$ to~$L$ is an
augmentation-preserving chain map taking~$K_q^*$ into~$L_q^*$ for
each~$q$. The resulting category of augmented directed complexes
is denoted $\mathbf{adc}$.

Given an augmented directed complex~$K$, we define an
$\omega$-category $\nu K$ functorially as follows. The members of
$\nu K$ are the double sequences
\[
 (x_0^-,x_0^+\mid x_1^-,x_1^+\mid\ldots\,)
\]
such that
\begin{align*}
 &x_q^\alpha\in K_q^*,\\
 &x_q^-=x_q^+=0\ \text{for all but finitely many values of $q$},\\
 &\epsilon x_0^-=\epsilon x_0^+=1,\\
 &x_q^+-x_q^-=\partial x_{q+1}^-=\partial x_{q+1}^+.
\end{align*}
The left and right identities $d_n^- x$ and $d_n^+ x$ of an
element
\[
 x=(x_0^-,x_0^+\mid x_1^-,x_1^+\mid\ldots\,)
\]
are given by
\[
 d_n^\alpha x=(x_0^-,x_0^+\mid\ldots\mid
 x_{n-1}^-,x_{n-1}^+\mid x_n^\alpha,x_n^\alpha\mid 0,0\mid\ldots\,).
\]
If $d_n^+x=d_n^- y$, say $d_n^+x=d_n^- y=z$, then the composite
$x\comp_n y$ is $x-z+y$.

We are especially interested in augmented directed complexes with
bases. These are augmented directed complexes of free abelian
groups with prescribed bases such that the distinguished
submonoids are generated, as monoids, by the prescribed basis
elements. It is convenient to work with the union of the bases for
the chain groups, a single graded set which we regard as the basis
for the entire augmented directed complex. In effect, we are
identifying a chain complex~$K$ with the direct sum
\[
 K_0\oplus K_1\oplus K_2\oplus\ldots.
\]
We note that an augmented directed complex can have at most one
basis, because a free abelian monoid has a unique basis
(consisting of its indecomposable non-zero elements). A globular
set~$X$ generates an augmented directed complex with basis equal
to the graded set~$X$ as follows: if $x\in X$ and $|x|=0$ then
$\epsilon x=1$; if $x\in X$ and $|x|>0$ then $\partial
x=\partial^+ x-\partial^- x$. We are interested in the cases when
$X$~is a simple globular set.

\begin{definition} \label{3.1}
A \emph{simple augmented directed complex} is an augmented
directed complex generated by a simple globular set.
\end{definition}

If $x$~is a positive-dimensional basis element in a simple
augmented directed complex, then $\partial^- x$ and $\partial^+ x$
are distinct basis elements. More generally let $K$ be any
augmented directed complex with a basis and let $x$ be any chain
in~$K$; then we will write
\[
 \partial x=\partial^+ x-\partial^- x,
\]
such that $\partial^- x$ and $\partial^+ x$ are sums of basis
elements with no common terms. This means that $\partial^+ x$ and
$\partial^- x$ are the `positive and negative parts of $\partial
x$'. The basis for~$K$ is called \emph{unital} if
\[
 \epsilon(\partial^-)^{|b|}b=\epsilon(\partial^+)^{|b|}b=1
\]
for every basis element~$b$. If the basis is unital, then we
associate an element~$\langle b\rangle$ of $\nu K$ called an
\emph{atom} to every basis element~$b$ as follows:
\[
 \langle b\rangle
 =\bigl(\,(\partial^-)^{|b|}b,\,(\partial^+)^{|b|}b\bigm|\ldots\bigm|
  \partial^- b,\,\partial^+ b\bigm|b,\,b\bigm|0,\,0\bigm|\ldots\,\bigl).
\]
In particular the basis of a simple augmented directed complex is
unital, because $(\partial^\alpha)^{|b|}b$ is a basis element for
every basis element~$b$ and because $0$-dimensional basis elements
have augmentation~$1$. If $b$~is a positive-dimensional basis
element in a simple augmented directed complex, then $\partial^-
b$ and $\partial^+ b$ are basis elements, so there are atoms
$\langle\partial^- b\rangle$ and $\langle\partial^+ b\rangle$. The
identities $\partial^\alpha\partial^-=\partial^\alpha\partial^+$
then give the following result.

\begin{proposition} \label{3.2}
If $b$~is an $n$-dimensional basis element in a simple augmented
directed complex with $n>0$, then
\[
 d_{n-1}^-\langle b\rangle=\langle\partial^- b\rangle,\
 d_{n-1}^+\langle b\rangle=\langle\partial^+ b\rangle.
\]
\end{proposition}

The main results of~\cite{SteOm} apply to augmented directed
complexes with bases when the bases are unital and are also
\emph{loop-free}, in the following sense: for $q\geq 0$ there is a
partial ordering~$<_q$ on the basis elements of degree at
least~$q$ such that $a<_q b$ if $a$~is a term in
$(\partial^-)^{|b|-q}b$ with $|b|>q$ or if $b$~is a term in
$(\partial^+)^{|a|-q}a$ with $|a|>q$. In practice we usually find
that the basis is \emph{strongly loop-free}; in this case there is
a partial ordering~$<_\N$ on the entire basis such that $a<_\N b$
if $a$~is a term in $\partial^- b$ or if $b$~is a term in
$\partial^+ a$. A strongly loop-free basis is loop-free, as the
terminology suggests, because one can get suitable partial
orderings~$<_q$ by restricting the partial ordering~$<_\N$. In
particular the basis of a simple augmented directed complex is
strongly loop-free, because the total ordering of the basis given
by its structure as a continuously graded ordered set has the
property required for~$<_\N$.

Let $K$ be an augmented directed complex with a loop-free unital
basis. The first main result of~\cite{SteOm} (Theorem~6.1) says
that the $\omega$-category $\nu K$ has a presentation of the
following kind: the generators are the atoms; for each
atom~$\langle b\rangle$ such that $|b|=n$ there are relations
$d_n^-\langle b\rangle=d_n^+\langle b\rangle=\langle b\rangle$;
for each atom~$\langle b\rangle$ such that $|b|=n>0$ there are
relations expressing $d_{n-1}^-\langle b\rangle$ and
$d_{n-1}^+\langle b\rangle$ as iterated composites of atoms (if
$d_{n-1}^\alpha\langle b\rangle$ can be expressed as an iterated
composite of atoms in more than one way, then we can choose any
such expression). When $K$~is a simple augmented directed complex,
we can use the formulae $d_{n-1}^\alpha\langle
b\rangle=\langle\partial^\alpha b\rangle$ of Proposition~\ref{3.2}
and we get the following result.

\begin{proposition} \label{3.3}
If $K$~is a simple augmented directed complex, then $\nu K$ is the
simple $\omega$-category generated by the basis for~$K$.
\end{proposition}

The second main result of~\cite{SteOm} (Theorem~5.11) says that
the functor~$\nu$ restricted to augmented directed complexes with
loop-free unital bases is a fully faithful embedding in the
category of $\omega$-categories. Since simple augmented directed
complexes have loop-free unital bases, a further restriction
produces the following result.

\begin{theorem} \label{3.4}
The restriction of~$\nu$ to the category of simple augmented
directed complexes is a fully faithful embedding with image
equivalent to the category~$\Theta$ of simple $\omega$-categories.
\end{theorem}

From this theorem we get a simple description of a category
equivalent to~$\Theta$. The objects are augmented chain complexes
of free abelian groups with prescribed non-empty finite ordered
bases $b_0,\ldots,b_p$ such that $|b_0|=|b_p|=0$ and such that
$|b_q|-|b_{q-1}|=\pm 1$ for $1\leq q\leq p$. The augmentation is
such that $\epsilon b_q=1$ for $|b_q|=0$. The boundary is such
that $\partial b_q=\partial^+ b_q-\partial^- b_q$ for $|b_q|>0$,
where $\partial^+ b_q$ is the last $(|b_q|-1)$-dimensional basis
element before~$b_q$ and where $\partial^- b_q$ is the first
$(|b_q|-1)$-dimensional basis element after~$b_q$. The morphisms
are the augmentation-preserving chain maps taking sums of
prescribed basis elements to sums of prescribed basis elements.

We will now show how to get a more combinatorial description of
these morphisms. Let $K$ be a simple augmented directed complex. A
finite, possibly empty, sequence of basis elements
$(a_1,\ldots,a_p)$ will be called \emph{ordered} if $a_1\leq
a_2\leq\ldots\leq a_p$. To an ordered sequence of $n$-dimensional
basis elements $(a_1,\ldots,a_p)$ we will associate the sum
$a_1+\ldots+a_p$; in this way the elements of~$K_n^*$ are
identified with the ordered sequences of $n$-dimensional basis
elements.

A sequence of $n$-dimensional basis elements $(a_1,\ldots,a_p)$
will be called \emph{separated} if it is ordered and if each
consecutive pair of terms encloses a basis element of lower
dimension; in particular a sequence of $0$-dimensional basis
elements is separated if and only if it is empty or a singleton.
Note that the members of a separated sequence are distinct.

The set of separated sequences of $n$-dimensional elements will be
partially ordered such that
$(a_1',\ldots,a_p')\leq(a_1'',\ldots,a_q'')$ if and only if the
following conditions hold:
\begin{align*}
 &p=q;\\
 &\text{$a_i'\leq a_i''$ for $1\leq i\leq p$};\\
 &\text{for $1\leq i\leq p$ the elements between $a_i'$~and~$a_i''$
 are at least $n$-dimensional.}
\end{align*}
Equivalently, for separated sequence of $n$-dimensional elements
$a'$~and~$a''$, one has $a'\leq a''$ if and only if $a''-a'$ is
the boundary of a member of~$K_{n+1}^*$.

Let
\[
 a'=(a_1',\ldots,a_p'),\quad a''=(a_1'',\ldots,a_p'')
\]
be separated sequences of $n$-dimensional elements such that
$a'\leq a''$. For $1\leq i\leq p$ let
$\{a_{i,0},\ldots,a_{i,q(i)}\}$ be the set of $n$-dimensional
elements~$a$ with $a_i'\leq a\leq a_i''$ indexed so that
\[
 a_i'=a_{i,0}<a_{i,1}<\ldots<a_{i,q(i)}=a_i''.
\]
A sequence of $(n+1)$-dimensional elements will be said to
\emph{bridge} $a'$~and~$a''$ if it has the form
\[
 (b_{1,1},\ldots,b_{1,q(1)},\ldots,b_{p,1},\ldots,b_{p,q(p)})
\]
such that $a_{i,j-1}<b_{i,j}<a_{i,j}$ for all $i$~and~$j$.
Equivalently, $b$~bridges $a'$~and~$a''$ if and only if $b$~is a
member of~$K_{n+1}^*$ such that $\partial b=a''-a'$. Note that a
sequence bridging two separated sequences is itself separated.

In these terms we get the following results.

\begin{theorem} \label{3.5}
Let $f\colon K\to L$ be a morphism of augmented directed complexes
be\-tween simple augmented directed complexes $K$~and~$L$, and let
$a$ be a basis element in~$K$. Then $fa$ is a separated sequence
of basis elements in~$L$.
\end{theorem}

\begin{proof}
The proof is by induction on~$|a|$.

If $|a|=0$ then $\epsilon fa=\epsilon a=1$, so $fa$ is a singleton
and is therefore separated.

If $|a|>0$ then $f\partial^- a$ and $f\partial^+ a$ are separated
by the inductive hypothesis and $fa$ bridges $f\partial^- a$ and
$f\partial^+ a$ because $\partial fa=f\partial
a=f\partial^+a-f\partial^- a$, so $fa$ is separated.

This completes the proof.
\end{proof}

\begin{theorem} \label{3.6}
Let $K$~and~$L$ be simple augmented directed complexes and let $f$
be a function assigning a separated sequence of $n$-dimensional
basis elements in~$L$ to every $n$-dimensional basis element
in~$K$. Then $f$~extends to a morphism of augmented directed
complexes if and only if the following conditions are satisfied.

(i) Let $a_1,\ldots,a_p$ be the zero-dimensional basis elements
in~$K$ indexed so that $a_1<\ldots<a_p$. Then the $fa_i$ are
singleton sequences such that $fa_1\leq\ldots\leq fa_p$.

(ii) For $n>0$ let $a'$~and~$a''$ be $(n-1)$-dimensional basis
elements in~$K$ such that $a'<a''$ and such that the elements
between $a'$~and~$a''$ have dimension at least~$n$. Let
$a_1,\ldots,a_p$ be the $n$-dimensional elements between
$a'$~and~$a''$ indexed so that $a'<a_1<\ldots<a_p<a''$. Then the
$fa_i$ are sequences bridging $fa'$ and $fa''$ such that
$fa_1\leq\ldots\leq fa_p$.
\end{theorem}

\begin{proof}
Condition~(i) is equivalent to saying that $f$~is
augmentation-preserving and that $f\partial K_1^*\subset\partial
L_1^*$. Condition~(ii) is equivalent to saying that $\partial
f=f\partial$ on~$K_n$ and that $f\partial K_{n+1}^*\subset\partial
L_{n+1}^*$. These conditions are clearly necessary and sufficient
for $f$~to give a morphism of simple augmented directed complexes.
\end{proof}

From Theorem~\ref{3.5}, because the members of separated sequences
are distinct, morphisms of augmented directed complexes between
simple augmented directed complexes can be regarded as functions
taking basis elements to sets of basis elements. Theorem~\ref{3.6}
therefore gives a combinatorial description of the morphisms
equivalent to the original definition~\cite{J} (see
also~\cite{MZ}). Theorem~\ref{3.6} also shows that the morphisms
can be constructed inductively: if one has constructed a morphism
in dimensions less than~$n$ satisfying the given conditions, then
one can always extend it to higher dimensions.

\section{The category of finite discs} \label{S4}

Makkai and Zawadowski \cite{MZ} have shown that Joyal's category
of finite discs \cite{J} is the opposite category to the category
of simple $\omega$-categories. Since morphisms between simple
$\omega$-categories can be represented by chain maps between
finitely generated free chain complexes, it follows that morphisms
between finite discs can be represented by cochain maps between
the dual cochain complexes. We will now give the details.

Recall that simple augmented directed complexes are augmented
chain complexes of finitely generated free abelian groups with
prescribed bases. Using dual bases, we see that the dual cochain
complexes are coaugmented cochain complexes of finitely generated
free abelian groups also with prescribed bases. It is clear that
the duals of the augmentation-preserving chain maps are the
coaugmentation-preserving cochain maps. Also, a chain map takes
sums of prescribed basis elements to sums of prescribed basis
elements if and only if its dual does the same for the dual bases
(in matrix terms, the entries of a matrix are nonnegative if and
only if the entries of its transpose are nonnegative). We
therefore get an easily described category equivalent to the
category of finite discs as follows.

The objects are coaugmented cochain complexes of free abelian
groups with prescribed non-empty finite bases $c_0,\ldots,c_p$
such that $|c_0|=|c_p|=0$ and such that $|c_q|-|c_{q-1}|=\pm 1$
for $1\leq q\leq p$. The coaugmentation~$\eta$ is such that
$\eta(1)$ is the sum of the $0$-dimensional basis elements. The
coboundary~$\delta$ is given by
\[
 \delta c_q=\delta^+ c_q-\delta^- c_q,
\]
where $\delta^+ c_q$ is the sum of the $(|c_q|+1)$-dimensional
elements~$c_r$ with $r<q$ such that $|c_t|>|c_q|$ for $r\leq t<q$,
and where $\delta^- c_q$ is the sum of the $(|c_q|+1)$-dimensional
elements~$c_s$ with $q<s$ such that $|c_t|>|c_q|$ for $q<t\leq s$.
The morphisms are the coaugmentation-preserving cochain maps
taking sums of prescribed basis elements to sums of prescribed
basis elements.

\section{Wreath products over the simplex category} \label{S5}

For $n=0,1,2,\ldots\,$, let $\Theta_n$ be the full subcategory
of~$\Theta$ consisting of the $n$-categories. We recall that an
$n$-category is an $\omega$-category in which every element is an
identity for~$\comp_n$; a simple $\omega$-category is therefore an
object of~$\Theta_n$ for which the dimensions in the associated
continuously graded ordered set do not exceed~$n$. The
categories~$\Theta_n$ form a chain
\[
 \Theta_0\subset\Theta_1\subset\Theta_2\subset\ldots
\]
and their union is the entire category~$\Theta$. Berger~\cite{Be}
has shown how to construct~$\Theta_n$ by using the functor
\[
 \mathcal{A}\mapsto\Delta\wr\mathcal{A}
\]
from categories to categories which takes each category to its
wreath product over the simplex category. He does this by
constructing equivalences $\Delta\wr\Theta_{n-1}\to\Theta_n$, from
which it follows that $\Theta_n$~is equivalent to the $n$-fold
iterated wreath product
\[
 \Delta\wr\ldots\wr\Delta.
\]
We will generalise his result by constructing a functor
\[
 V\colon\Delta\wr\mathbf{adc}\to\mathbf{adc}
\]
and showing that its restriction to an appropriate subcategory is
fully faithful.

We begin by recalling the definition of the wreath product
category $\Delta\wr\mathcal{A}$ for an arbitrary
category~$\mathcal{A}$. The objects of $\Delta\wr\mathcal{A}$ are
the pairs $(m,A)$, where $m$~is a nonnegative integer and where
$A=(A^1,\ldots,A^m)$ is an ordered $m$-tuple of objects
of~$\mathcal{A}$. The morphisms in $\Delta\wr\mathcal{A}$ from
$(m,A)$ to $(n,B)$ are the pairs $(\phi,f)$, where
$\phi=\bigl(\phi(0),\ldots,\phi(m)\bigr)$ is an ordered
$(m+1)$-tuple of integers with
\[
 0\leq\phi(0)\leq\phi(1)\leq\ldots\leq\phi(m)\leq n
\]
and where
\[
 f=(f_1^{\phi(0)+1},\ldots,f_1^{\phi(1)},
 f_2^{\phi(1)+1},\ldots,f_2^{\phi(2)},\ldots,
 f_m^{\phi(m-1)+1},\ldots,f_m^{\phi(m)})
\]
is an ordered $[\phi(m)-\phi(0)]$-tuple of morphisms $f_i^j\colon
A^i\to B^j$ in~$\mathcal{A}$. Composition in
$\Delta\wr\mathcal{A}$ is given by
\[
 (\psi,g)\circ(\phi,f)=(\psi\circ\phi,g\circ f),
\]
where $\psi\circ\phi(i)=\psi\bigl(\phi(i)\bigr)$ and where
$(g\circ f)_i^k=g_j^k\circ f_i^j$ for
\[
 \psi\circ\phi(i-1)\leq\psi(j-1)\leq k-1<k\leq\psi(j)
 \leq\psi\circ\phi(i).
\]
If $\mathcal{A}$ is the category with one object and one morphism,
then one sees that $\Delta\wr\mathcal{A}$ is the simplex
category~$\Delta$ itself. It is clear that the wreath product
construction $\Delta\wr-$ is functorial.

We will now construct the functor
\[
 V\colon\Delta\wr\mathbf{adc}\to\mathbf{adc}.
\]
As before, it is convenient to regard an augmented directed
complex~$K$ as the direct sum of its chain groups~$K_q$, so that
\[
 K=K_0\oplus K_1\oplus\ldots.
\]

Let $(m,K)$ be an object of $\Delta\wr\mathbf{adc}$ with
$K=(K^1,\ldots,K^m)$. For $0\leq i\leq m$, let $\Z p^i$ be a free
abelian group with a single basis element~$p^i$. For $1\leq i\leq
m$ let $sK^i$ be an abelian group isomorphic to~$K^i$ and let
$s\colon K^i\to sK^i$ be an isomorphism. Then
\[
 V(m,K)
 =\Z p^0\oplus sK^1\oplus\Z p^1\oplus\ldots
 \oplus\Z p^{m-1}\oplus sK^m\oplus\Z p^m
\]
with the following structure. The degrees are given by $|p^i|=0$
and $|s x|=|x|+1$. The boundary is given as follows: if $x\in K^i$
with $|x|>0$ then $\partial sx=s\partial x$; if $x\in K^i$ with
$|x|=0$ then
\[
 \partial sx=(\epsilon x)(p^i-p^{i-1});
\]
and finally $\partial p^i=0$. The augmentation is given by
$\epsilon p^i=1$. The distinguished submonoid of $V(m,K)$ is
generated by the elements~$p^i$ and by the images under~$s$ of the
distinguished submonoids of the~$K^i$.

Now let $(\phi,f)\colon(m,K)\to(n,L)$ be a morphism in
$\Delta\wr\mathbf{adc}$. Then
\[
 V(\phi,f)\colon V(m,K)\to V(n,L)
\]
is the homomorphism given as follows: if $0\leq i\leq m$ then
$p^i\mapsto p^{\phi(i)}$; if $x\in K^i$ then
\[
 sx\mapsto\sum_{\phi(i-1)<j\leq\phi(i)}sf_i^j x.
\]
It is straightforward to check that $V(\phi,f)$ is a morphism of
augmented directed complexes, and that $V$~is a functor.

\begin{example} \label{5.1}
Let $K^1,\ldots,K^m$ be simple augmented directed complexes and
let the sequences of dimensions for the underlying continuously
graded ordered sets be $t^1,\ldots,t^m$. Then $V(m,K)$ is a simple
augmented directed complex such that the sequence of dimensions
for the underlying continuously graded ordered set is
\[
 (0,st^1,0,st^2,0,\ldots,0,st^k,0),
\]
where the sequence $st^i$ is obtained from the sequence~$t^i$ by
adding~$1$ to each term.
\end{example}

\begin{example} \label{5.2}
Let $K^1,\ldots,K^m$ be the cellular chain complexes of cell
complexes~$X^i$ with prescribed orientations for the cells, made
into augmented directed complexes by taking the oriented cells as
bases. Let $I$ be the closed interval $[0,1]$, let
$P^0,\ldots,P^m$ be one-point spaces, and let $V(m,X)$ be obtained
from the disjoint union of the spaces
\[
 P^0,\ X^1\times I,\ P^1,\ \ldots,\ P^{m-1},\ X^m\times I,\ P^m
\]
by identifying $X^i\times\{0\}$ with~$P^{i-1}$ and
$X^i\times\{1\}$ with~$P^i$. Thus $V(m,X)$ is a chain got by
joining the unreduced suspensions of the spaces~$X^i$ together;
alternatively, $V(m,X)$ is the homotopy colimit of the unique
diagram of the form
\[
 P^0\leftarrow X^1\rightarrow P^1\leftarrow\ldots
 \rightarrow P^{m-1}\leftarrow X^m\rightarrow P^m.
\]
There is an obvious cell structure on $V(m,X)$: the cells are the
products $c\times I$ for $c$ a cell in some~$X^i$ together with
the points~$P^i$. If one orientates the cells correctly, then
$V(m,K)$ becomes the cellular chain complex of $V(m,X)$.
\end{example}

The functor~$V$ behaves well on bases.

\begin{proposition} \label{5.3}
Let $(m,K)$ be an object of $\Delta\wr\mathbf{adc}$ such that the
augmented directed complexes~$K^i$ have bases. Then $V(m,K)$ has a
basis. If the bases for the~$K^i$ are unital then the basis for
$V(m,K)$ is unital. If the bases for the~$K^i$ are loop-free then
the basis for $V(m,K)$ is loop-free. If the bases for the~$K^i$
are strongly loop-free then the basis for $V(m,K)$ is strongly
loop-free.
\end{proposition}

\begin{proof}
It is straightforward to check that $V(m,K)$ has a basis
consisting of the zero-dimensional elements~$p^i$ together with
the images under~$s$ of the bases for the terms~$K^i$.

Suppose that the bases for the~$K^i$ are unital. By construction,
$\epsilon p^i=1$. If $b$~is a basis element for~$K^i$, then
\[
 \epsilon(\partial^-)^{|sb|}sb
 =\epsilon\partial^-(\partial^-)^{|b|}sb
 =\epsilon\partial s(\partial^-)^{|b|}b
 =\epsilon[\epsilon(\partial^-)^{|b|}b]p^{i-1}
 =\epsilon(1p^{i-1})
 =1
\]
and $\epsilon(\partial^+)^{|sb|}sb=\epsilon p^i=1$ similarly.
Therefore the basis for $V(m,K)$ is unital.

Suppose that the bases for the~$K^i$ are loop-free, with partial
orderings~$<_q$. For $q>0$ one gets a partial ordering~$<_q$ in
$V(m,K)$ as required for loop-freeness by taking the union of the
images of the partial orderings~$<_{q-1}$ under the isomorphisms
$s\colon K^i\to sK^i$. One also gets a partial ordering~$<_0$ in
$V(m,K)$ with the required property as follows:
\begin{align*}
 &\text{$p^i<_0 p^j$ for $i<j$,}\\
 &\text{$p^i<_0 sb$ for $b\in K^j$ with $i<j$,}\\
 &\text{$sa<_0 p^j$ for $a\in K^i$ with $i\leq j$,}\\
 &\text{$sa<_0 sb$ for $a\in K^i$ and $b\in K^j$ with $i<j$.}
\end{align*}
Therefore the basis for $V(m,K)$ is loop-free.

Now suppose that the bases for the~$K^i$ are strongly loop-free
with partial orderings~$<_\N$. Then the basis for $V(m,K)$ is
strongly loop-free under the partial ordering~$<_\N$ given by
\begin{align*}
 &\text{$p^i<_\N p^j$ for $i<j$,}\\
 &\text{$p^i<_\N sb$ for $b\in K^j$ with $i<j$,}\\
 &\text{$sa<_\N p^j$ for $a\in K^i$ with $i\leq j$,}\\
 &\text{$sa<_\N sb$ for $a\in K^i$ and $b\in K^j$ with
 $i<j$,}\\
 &\text{$sa<_\N sb$ for $a<_\N b$ in some~$K^i$.}
\end{align*}

This completes the proof.
\end{proof}

Let $\Phi$ be the full subcategory of $\mathbf{adc}$ consisting of
non-zero augmented directed chain complexes with loop-free unital
bases. Because of Proposition~\ref{5.3}, the functor $V$ maps
$\Delta\wr\Phi$ into~$\Phi$, and the main result on the
functor~$V$ says that the restriction
\[
 V\colon\Delta\wr\Phi\to\Phi
\]
is fully faithful. We need two subsidiary results.

\begin{proposition} \label{5.4}
Let $K$ be a non-zero augmented directed complex with a unital
basis. Then $K$~has at least one zero-dimensional basis element,
and every zero-dimensional basis element has augmentation~$1$.
\end{proposition}

\begin{proof}
Since $K$~is non-zero, there is at least one basis element~$b$.
For this element we have $\epsilon(\partial^-)^{|b|}b=1$, so
$(\partial^-)^{|b|}b$ is a non-zero chain of dimension zero. It
follows that there is at least one zero-dimensional basis element.
Finally, if $a$~is a zero-dimensional basis element then $\epsilon
a=\epsilon(\partial^-)^{|a|}a=1$.
\end{proof}

\begin{proposition} \label{5.5}
Let $K$~and~$L$ be augmented directed complexes with bases such
that the basis for~$L$ is loop-free and unital. If $f\colon K\to
L$ is a chain map taking sums of basis elements to sums of basis
elements such that $\epsilon\circ f=0$, then $f=0$.
\end{proposition}

\begin{proof}
It suffices to show that $fa=0$ for every basis element~$a$
in~$K$, and we will use induction on~$|a|$.

Suppose that $|a|=0$. Then
\[
 fa=b_1+\ldots+b_k
\]
for some basis elements~$b_i$ in~$K$. We have $\epsilon b_i=1$ for
each~$i$, because the basis for~$L$ is unital, so $\epsilon fa=k$.
But $\epsilon fa=0$, so $k=0$, which means that $fa=0$.

Now suppose that $|a|=q+1>0$. Again we have
\[
 fa=b_1+\ldots+b_k
\]
for some basis elements~$b_i$, and we must show that $k=0$.
Suppose therefore that $k>0$. Since the basis for~$L$ is
loop-free, we can choose~$b_i$ such that $b_j<_q b_i$ is not true
for any~$j$. Since $\epsilon(\partial^-)^{q+1}b_i=1$, we must have
$\partial^- b_i\neq 0$. But a basis element which is a term in
$\partial^- b_i$ cannot be cancelled by a term in $\partial^+ b_j$
for any~$j$, because we do not have $b_j<_q b_i$, so we get
\[
 f\partial a=\partial fa=\partial b_1+\ldots+\partial b_k\neq 0,
\]
contrary to the inductive hypothesis. Therefore $fa=0$.

This completes the proof.
\end{proof}

The main result is now as follows.

\begin{theorem} \label{5.6}
Let $\Phi$ be the full subcategory of $\mathbf{adc}$ consisting of
the non-zero augmented directed complexes with loop-free unital
bases. Then the functor
\[
 V\colon\Delta\wr\Phi\to\Phi
\]
is fully faithful.
\end{theorem}

\begin{proof}
Let $(m,K)$ and $(n,L)$ be objects of $\Delta\wr\Phi$, and let
\[
 F\colon V(m,K)\to V(n,L)
\]
be a morphism in~$\Phi$. We must show that $F=V(\phi,f)$ for a
unique morphism $(\phi,f)\colon(m,K)\to(n,L)$ in $\Delta\wr\Phi$.

As an abelian group, $V(m,K)$ is generated by the elements~$p^i$
and the subgroups $sK^i$, so $F$~is determined by its values
on~$p^i$ and on $sK^i$. Since $F$~is an augmentation-preserving
chain map taking sums of basis elements to sums of basis elements,
$F p^i$ must be a sum of zero-dimensional basis elements such that
$\epsilon F p^i=1$. This forces $F p^i$ to be a single basis
element, so
\[
 F p^i=p^{\phi(i)}
\]
for some $\phi(i)$ with $0\leq\phi(i)\leq n$. Also, since $F$~is a
chain map taking sums of basis elements to sums of basis elements,
we must have
\[
 Fsx=\sum_{j=1}^n sf_i^j x\ \text{for $x\in K^i$},
\]
where the $f_i^j\colon K^i\to L^j$ are uniquely determined chain
maps taking sums of basis elements to sums of basis elements. To
show that $F=V(\phi,f)$ for a unique morphism
$(\phi,f)\colon(m,K)\to(n,L)$ in $\Delta\wr\Phi$ it now suffices
to show that
\[
 \phi(0)\leq\phi(1)\leq\ldots\leq\phi(m),
\]
that $f_i^j$~is augmentation-preserving for
$\phi(i-1)<j\leq\phi(i)$, and that $f_i^j=0$ otherwise.

We first show that $\phi(i-1)\leq\phi(i)$, using the assumption
that $K^i\neq 0$. By Proposition~\ref{5.4} there is a
zero-dimensional basis element~$a$ in~$K^i$, and $\epsilon a=1$.
It follows that
\[
 F\partial sa=F(\epsilon a)(p^i-p^{i-1})=p^{\phi(i)}-p^{\phi(i-1)}
\]
and that
\[
 \partial Fsa
 =\partial\sum_{j=1}^n sf_i^j a
 =\sum_{j=1}^n(\epsilon f_i^j a)(p^j-p^{j-1}).
\]
Since $f_i^j a$ is a sum of basis elements, we must have $\epsilon
f_i^j a\geq 0$ for all~$j$. Since $F\partial sa=\partial Fsa$, we
must have $\phi(i-1)\leq\phi(i)$.

Next we compute $\epsilon\circ f_i^j$. Let $x$ be a
zero-dimensional chain in~$K^i$. Equating $F\partial sx$ and
$\partial Fsx$ just as in the previous paragraph, we get
\[
 (\epsilon x)(p^{\phi(i)}-p^{\phi(i-1)})
 =\sum_{j=1}^n(\epsilon f_i^j x)(p^j-p^{j-1}).
\]
For $\phi(i-1)<j\leq\phi(i)$ we get $\epsilon f_i^j x=\epsilon x$,
so that $f_i^j$~is augmentation-preserving, and for other values
of~$j$ we get $\epsilon f_i^j x=0$, so that $f_i^j=0$ by
Proposition~\ref{5.5}.

This completes the proof.
\end{proof}

Let us now consider the categories~$\Theta_n$ in the filtration
of~$\Theta$ given by
\[
 \Theta_0\subset\Theta_1\subset\Theta_2\subset\ldots.
\]
The category~$\Theta$ is equivalent to a full subcategory of
$\mathbf{adc}$, and its objects are non-zero augmented directed
complexes with loop-free unital bases. It follows from
Theorem~\ref{5.6} that $V$~induces a fully faithful functor from
$\Delta\wr\Theta$ to strict $\omega$-categories. From
Example~\ref{5.1} we see that this functor restricts to an
equivalence $\Delta\wr\Theta_{n-1}\cong\Theta_n$; this equivalence
is also visible in Theorem~\ref{3.6}. Now an object of~$\Theta_0$
is an augmented directed complex with a basis consisting of a
single zero-dimensional element of augmentation~$1$. Clearly, if
$K$~and~$L$ are objects of~$\Theta_0$ then there is a unique
morphism of augmented directed complexes from~$K$ to~$L$. It
follows that $\Theta_0$~is equivalent to the category with a
single object and a single morphism. For $n>0$ it then follows
that $\Theta_n$ is equivalent to the $n$-fold iterated wreath
product $\Delta\wr\ldots\wr\Delta$; we have therefore recovered
Berger's result (\cite{Be}, Theorem~3.7). The inclusion functor
$\Theta_0\to\Theta_1$ is equivalent to the functor
$\Theta_0\to\Delta\wr\Theta_0$ sending the objects of~$\Theta_0$
to the unique object of the form $(0,A)$ in $\Delta\wr\Theta_0$,
and the inclusions $\Theta_{n-1}\to\Theta_n$ are obtained from
this by repeated application of the wreath product functor. Up to
equivalence, the categories~$\Theta_n$ and the inclusion functors
$\Theta_{n-1}\to\Theta_n$ can therefore be expressed entirely in
terms of the simplex category and wreath products. The entire
category~$\Theta$ can then be expressed in this way as well, since
it is the colimit of the sequence
\[
 \Theta_0\to\Theta_1\to\Theta_2\to\ldots.
\]

\begin{remark} \label{5.7}
In particular we have an isomorphism $\Delta\cong\Theta_1$ giving
a fully faithful embedding of the simplex category~$\Delta$ in the
category of augmented directed complexes $\mathbf{adc}$. There is
a completely different embedding of~$\Delta$ in $\mathbf{adc}$
sending the objects of~$\Delta$ to the simplicial chain complexes
of the standard simplexes (see~\cite{SteOr}). This embedding is
not full. The corresponding $\omega$-categories are Street's
orientals (see~\cite{StrAlg}). The oriental corresponding to the
$n$-dimensional simplex is a simple $\omega$-category only for
$n=0$ and for $n=1$.
\end{remark}


\begin{thebibliography}{1}

\bibitem{Ba}
M.~A. Batanin, \emph{Monoidal globular categories as a natural
environment for
  the theory of weak {$n$}-categories}, Adv. Math. \textbf{136} (1998), no.~1,
  39--103.

\bibitem{BS}
Michael Batanin and Ross Street, \emph{The universal property of
the
  multitude of trees}, J. Pure Appl. Algebra \textbf{154} (2000), no.~1-3,
  3--13.

\bibitem{Be}
Clemens Berger, \emph{Iterated wreath product of the simplex
category and
  iterated loop spaces}, Adv. Math., to appear.

\bibitem{J}
Andr\'e Joyal, \emph{Disks, duality and $\theta$-categories},
unpublished, 1997.

\bibitem{MZ}
Mihaly Makkai and Marek Zawadowski, \emph{Duality for simple
  {$\omega$}-categories and disks}, Theory Appl. Categ. \textbf{8} (2001),
  114--243 (electronic).

\bibitem{SteOm}
Richard Steiner, \emph{Omega-categories and chain complexes},
Homology Homotopy
  Appl. \textbf{6} (2004), no.~1, 175--200 (electronic).

\bibitem{SteOr}
\bysame, \emph{Orientals}, arXiv:math.CT/0601383, 2006.

\bibitem{StrAlg}
Ross Street, \emph{The algebra of oriented simplexes}, J. Pure
Appl. Algebra
  \textbf{49} (1987), no.~3, 283--335.

\bibitem{StrPet}
\bysame, \emph{The petit topos of globular sets}, J. Pure Appl.
Algebra
  \textbf{154} (2000), no.~1-3, 299--315.

\end{thebibliography}

\providecommand{\bysame}{\leavevmode\hbox
to3em{\hrulefill}\thinspace}
\providecommand{\MR}{\relax\ifhmode\unskip\space\fi MR }
\providecommand{\MRhref}[2]{%
  \href{http://www.ams.org/mathscinet-getitem?mr=#1}{#2}
} \providecommand{\href}[2]{#2}

\end{document}